\documentclass[11pt]{article}
\newcommand{\documentdate}{28 I 2025}

\usepackage{a4wide}
\usepackage{graphicx}
\usepackage{subfigure}
\usepackage{amsmath}
\usepackage{amssymb}
\usepackage{mathtools,xcolor}
\usepackage{amsthm}
\usepackage{dsfont}

\newcommand{\beqn}[1]{\begin{equation}\label{#1}}
\newcommand{\eeqn}{\end{equation}}
\newcommand{\req}[1]{(\ref{#1})}
\newcommand{\ms}{\;\;\;\;}
\newcommand{\tim}[1]{\;\; \mbox{#1} \;\;}

\newcommand{\bpr}{{\bf Proof.} \hspace{1.5mm}}
\newcommand{\epr}{\hfill $\Box$ \vspace*{1em}}

\newcommand{\ii}[1]{\{ 1, \ldots, #1 \}}
\newcommand{\iiz}[1]{\{ 0, \ldots, #1 \}}
\newcommand{\iibe}[2]{\{ #1, \ldots, #2 \}}
\newcommand{\gamlow}{\gamma_{\rm low}}
\newcommand{\calA}{{\cal A}} 
\newcommand{\calE}{{\cal E}} 
\newcommand{\calF}{{\cal F}}

\newcommand{\calO}{{\cal O}} 
\renewcommand{\Re}{\hbox{I\hskip -2pt R}}
\newcommand{\smallRe}{\hbox{\footnotesize I\hskip -2pt R}}
\newcommand{\bigfrac}[2]{\frac{\displaystyle #1}{\displaystyle #2}}
\newcommand{\bigsum}{\displaystyle \sum}

\newcommand{\sfrac}[2]{{\scriptstyle \frac{#1}{#2}}}
\newcommand{\kap}[1]{\kappa_{\mbox{\rm\tiny #1}}}
\newcommand{\eqdef}{\stackrel{\rm def}{=}}
\newcommand{\al}[1]{{\footnotesize{\sf #1}}}
\newcommand{\tal}[1]{{\normalsize {\sf #1}}}
\newcommand{\half}{\sfrac{1}{2}}

\newcommand{\quarter}{\sfrac{1}{4}}

\newcommand{\Flow}{F_{\rm low}}

\newcommand{\E}[1]{\mathbb{E}\!\left[#1 \right]}

\newcommand{\Econd}[2]{\mathbb{E}_{#2}\!\left[{#1}\right]}
\newcommand{\Pcond}[2]{\mathbb{P}_{#2}\!\left[#1\right]}
\newcommand{\Vcond}[2]{\mathrm{Var}_{#2}\left[{#1}\right]}

\DeclareMathOperator*{\average}{average}
\newcommand{\comment}[1]{}
\newcommand{\avggik}[2]{\frac{1}{#1}\sum_{j=0}^{#2} |g_{i,j}|}
\newcommand{\indic}[1]{\mathds{1}_{#1}}

\theoremstyle{plain}
\newtheorem{theorem}{Theorem}[section]

\newtheorem{lemma}[theorem]{Lemma}

\newcommand{\numsection}[1]{\section{#1}\setcounter{equation}{0}}
\newcommand{\appnumsection}[1]{\section*{#1}\setcounter{equation}{0}
  \renewcommand{\theequation}{A.\arabic{equation}}
  \renewcommand{\thetheorem}{A.\arabic{theorem}}
  \renewcommand{\thetable}{A.\arabic{table}}
  \renewcommand{\thefigure}{A.\arabic{figure}}
  \renewcommand{\thesection}{A} }
\renewcommand{\theequation}{\arabic{section}.\arabic{equation}}
\renewcommand{\thetable}{\arabic{section}.\arabic{table}}
\renewcommand{\thefigure}{\arabic{section}.\arabic{figure}}

\theoremstyle{definition}

\newtheorem{assumption}[theorem]{Assumption}
\theoremstyle{remark}

\newcounter{algo}[section]
\renewcommand{\thealgo}{\thesection.\arabic{algo}}
\newcommand{\algo}[3]{\refstepcounter{algo}
\begin{center}\begin{figure}[htbp]
\framebox[\textwidth]{
\parbox{0.95\textwidth} {\vspace{\topsep}
{\bf Algorithm \thealgo : #2}\label{#1}\\
\vspace*{-\topsep} \mbox{ }\\
{#3} \vspace{\topsep} }}
\end{figure}\end{center}}

\title{Complexity and performance  for
       two classes\\ of noise-tolerant first-order algorithms}

\author{
   S. Gratton%
   \thanks{Universit\'e de Toulouse, INP, IRIT, Toulouse, France. Email:
     serge.gratton@enseeiht.fr. Work partially supported by 3IA Artificial and
     Natural Intelligence Toulouse Institute (ANITI), French "Investing for the Future
     - PIA3" program under the Grant agreement ANR-19-PI3A-0004"}, 
   ~S. Jerad%
   \thanks{ANITI, Universit\'e de Toulouse, INP, IRIT, Toulouse, France. Email:
     sadok.jerad@enseeiht.fr}
   ~and Ph. L. Toint%
   \thanks{NAXYS, University of Namur, Namur, Belgium. Email:
     philippe.toint@unamur.be.
     Partially supported by ANITI.}
}

\date{\documentdate}

\begin{document}

\maketitle

\begin{abstract}
Two classes of algorithms for optimization in the presence of noise are presented, that
do not require the evaluation of the objective function. The first generalizes the well-
known Adagrad method. Its complexity is then analyzed as a function of its parameters.
A second class of algorithms is then derived whose complexity
is at least as good as that of the first class. Initial numerical experiments on finite-sum
problems arising from deep-learning applications suggest  that methods of the second
class may outperform those of the first.
\end{abstract}


\numsection{Introduction}
\label{introduction}

Minimization algorithms which can handle noisy evaluations of the objective
function and/or gradients have generated a significant amount of research in
the last few years
\cite{toint2020stochastic,Bellavia2021,BellGuriMoriToin21b,
  Bellavia2022,Berahas2021,Blanchet2019,
  CartisScheinberg,Chen2017,defossez2022a,duchi11adagrad,KingBa15,
  LiOrabona2018,Paquette2020,ReddiKK18,TielHint12,
  WardWuBott19,ZhouTangYangCaoGu24,ZouShenJieSunLiu19a}.
Interestingly, a number of these contributions
\cite{toint2020stochastic,Bellavia2021,BellGuriMoriToin21b,
  Bellavia2022,Berahas2021,Blanchet2019,
  CartisScheinberg,Chen2017,Paquette2020}
indicate that, when the (noisy) objective function is evaluated,
its accuracy is significantly more critical to ensure convergence than that of
the computed (noisy) derivatives. This may be the reason why methods where the problem is avoided by
\emph{not} evaluating the objective function  
 (such as Adagrad
\cite{DuchHazaSing11}, RMSProp \cite{TielHint12}, Adam \cite{KingBa15} or
AMSGrad \cite{ReddKaleKuma18}), have become very popular in the context of
finite-sum minimization, where noise in the evaluation arises from sampling
among a very large number of terms. That such methods can be provably
convergent to first-order stationary points is quite remarkable,  and
the literature covering their theory is extensive. We now briefly
survey some of the contributions most relevant in our context.

\subsection{Related work}

Several authors have been able to prove global convergence rates, including
the recent contributions of
\cite{defossez2022a,Fawetal22,kavis2022,GadatGavara22,Wangetal23,
ZhouTangYangCaoGu24,LiOrabona2018,Fawetal23,AttiaKoren23,WardWuBott19},
where a global convergence analysis of the Adagrad method has been conducted
under different assumptions. The paper \cite{LiOrabona2018} provides
an analysis of a delayed Adagrad method which does not consider the
sampled gradient at iteration $k$ to compute its step size, and
further variants that are arbitrarily close to Adagrad and
require specific choices of hyperparameters. The first parameter-free
analysis of Adagrad was proposed in \cite{WardWuBott19} and later
revisited by \cite{defossez2022a}, where the dependence on some
parameters was improved. In both cases the bounds were given in
expectation. In \cite{kavis2022}, high-probability bounds are derived
for a class of methods that includes Adagrad. In
\cite{ZhouTangYangCaoGu24}, complexity analyses of a large class of
adaptive gradient methods (including Adagrad) are proposed, and
improved convergence rates are proved under the ``gradient sparsity''
assumption of the gradient iteration sequence. Note that all of 
\cite{WardWuBott19,defossez2022a,kavis2022,ZhouTangYangCaoGu24}
require  the sampled gradient to be bounded. This
requirement was circumvented in
\cite{Fawetal22,Wangetal23,Fawetal23,AttiaKoren23} by allowing
unbounded gradients boundedness using a new Lipschitz smoothness
condition \cite{Zhang2020Why} and different noise assumptions, see
\cite{Fawetal22,Wangetal23,Fawetal23,AttiaKoren23} for
more details.  

\subsection{Our contributions}

The present paper remains in the context of bounded gradients and extends some results of
\cite{WardWuBott19,defossez2022a} to achieve several goals.
\begin{enumerate}
\item The global rate of convergence result of \cite{defossez2022a} is
  shown to hold for an extended class of methods including the Adagrad algorithm.
\item Using the new analysis tools, a new class of methods is then
  proposed, whose global rate of convergence is shown to be very close
  to that of methods using (exact) function evaluations.
\item Numerical experiments with finite-sum problems arising from
  deep-learning applications indicate that method of the latter class
  may sometimes perform better than those of the former.
\end{enumerate}

The presentation is organized as follows. A general framework of first-order
trust-region algorithms is introduced in Section~\ref{section:algo}, in which two
classes of algorithms (one of them containing the Adagrad method) are defined and
analyzed (complexity-wise) in Sections~\ref{section:adag} and
\ref{section:divseries}, respectively. Numerical experiments in the
finite-sum minimization context are presented in Section~\ref{section:numerics}.
Some conclusions are finally outlined in Section~\ref{section:conclusion}.

\numsection{A first-order framework for minimizing noisy functions}
\label{section:algo}

We are interested in (approximately) solving the problem
\beqn{problem}
\min_{x\in\smallRe^n} F(x)
\eeqn
where $F$ is a function from $\Re^n$ to $\Re$ contaminated by noise.
Moreover, we assume that evaluating $F$ at any given $x$ to sufficient
accuracy is either impossible or too costly.  Evaluating a noisy gradient is
however possible\ldots and our only source of information about the problem.
While access to $F$ or its exact gradient is impossible, we nevertheless make
the following assumptions.

\begin{assumption}\label{AS.1}
The objective function $F(x)$ is continuously differentiable.
\end{assumption}
\begin{assumption}\label{AS.2}
Its exact gradient $G(x) \eqdef \nabla_x^1f(x)$ is Lipschitz
continuous with Lipschitz constant $L$, that is
\[
\|G(x)-G(y)\| \le L \|x-y\|
\]
for all $x,y\in \Re^n$.
\end{assumption}
\begin{assumption}\label{AS.3}
 There exists a constant $\Flow$ such that, for all $x$, $F(x)\ge \Flow$.
\end{assumption}

\noindent
A standard consequence of Assumption~\ref{AS.2} is that, for, any $x,s \in \Re^n$,
\beqn{lipsch-bound}
F(x+s) \leq F(x) + G(x)^Ts + \frac{L}{2} \|s\|^2
\eeqn
(see Lemma~2.1 in \cite{Bellavia2019} or Theorem~A.8.3 in \cite{Cartis2022-wb}, for instance).

\noindent
We now present a first-order \emph{adaptively scaled gradient} algorithmic framework
(\al{ASGRAD}), where, at iteration $k$, a \emph{noisy} gradient $g_k= g(x_k)$
is evaluated and a step $s_k$ defined that decreases the associated local
linear model and whose size is determined by componentwise ``scaling factors''
$w_{i,k}$ to be chosen at each iteration. Our framework is formally described
as follows. 

\algo{ASGRAD}{The \tal{ASGRAD} framework}{
\vspace*{-4mm}
\begin{description}
\item[Step 0: Initialization. ]
  $x_0$ and a constant $\gamma_{\rm
  low}\in (0,1]$  are given. Set $k=0$.
\item[Step 1: Step computation. ]
  Evaluate $g_k$ and set
  \beqn{sQ-def-a}
  s_k = \gamma_{k} s_k^L,
  \eeqn
  with
  \beqn{sL-def-a}
  s^L_{i,k} = -\frac{g_{i,k}}{w_{i,k}}
  \eeqn
  for a stepsize $\gamma_k \in [\gamlow, 1]$ and positive scaling factors $w_{i,k}$.

\item[Step 2: New iterate.]
  Define
  \beqn{xupdate-a}
  x_{k+1} = x_k + s_k,
  \eeqn
  increment $k$ by one and return to Step~1.
\end{description}
}
  
\noindent
We stress that $g_k$ (as evaluated in Step~1) is a noisy random gradient
evaluation. The algorithms of the \al{ASGRAD} framework therefore generate a stochastic process
\[
\{ x_k, g_k, \gamma_k, s^L_k, s_k \}
\]
on some probability space $(\Omega,\calF,\mathbb{P})$.  The associated expectation
operator will be denoted by $\E{\cdot}$ and $\Econd{\cdot}{k}$ will stand for the
conditional expectation knowing $\{g_0,\ldots,g_{k-1}\}$. All
algorithms in our framework may clearly be interpreted as variants of
Stochastic Gradient Descent, allowing for a variety of stepsize
(learning rate) rules.

We will, in what follows, assume that the noisy gradient $g_k$ is a bounded
non-biased estimator of the true gradient, that is
\noindent
\begin{assumption}\label {AS.4}
We have that, for all $k\geq0$, $\Econd{g_k}{k}
  = G(x_k)$. Moreover, there exists a constant $\kappa_g \geq 1$ such that
  $\|g_k\|_\infty \leq \kappa_g$ for all $k\geq 0$ and all realizations of the
  algorithm.
In addition, we assume that $\gamma_{k}$ is mesurable with respect to $\{g_0,\ldots,g_{k-1}\}$.
\end{assumption}

\noindent
This assumption that gradients are bounded maybe quite realistic in practice\footnote{Not to
mention that an infinite gradient is likely to crash the algorithm on
many machines.}, for instance when the iterates remain in a compact
subset of $\Re^n$ and has been extensively used in the analysis of 
stochastic first-order methods (see
\cite{WardWuBott19,defossez2022a,ZhouTangYangCaoGu24,Zaheeretal18,kavis2022})
immediately implies that 
\beqn{Gbounded} 
\|G(x_k)\|_\infty \leq \kappa_g \tim{ for all } k\geq 0.
\eeqn

\noindent
The assumption on the stepsize $\gamma_{k}$  is consistent with
current best practices in training deep neural networks: to achieve
state-of-the-art performance, the step size $\gamma_{k}$ is often
either set to a specific value  at the beginning and divided by a
fixed constant at (approximately) regular intervals
\cite{ZagoruykoK16} or periodically warm-restarted
\cite{LoschHutter17}.

\noindent
The reader has undoubtedly noted that we have not been very specific regarding
how the scaling factors $w_{i,k}$ are selected, and a whole range of options
is possible. This justifies our choice to consider \al{ASGRAD} as an
\emph{algorithmic framework}, covering many possible such choices. The rest of this paper is
devoted to the analysis of two specific classes of interest.

\numsection{An Adagrad-inspired class of \tal{ASGRAD} algorithms}\label{section:adag}

In the first considered \al{ASGRAD} class, the scaling factors are inspired by the
definition of the Adagrad algorithm \cite{DuchHazaSing11}.  More
specifically, we make the following additional assumptions.

\noindent
\begin{assumption}\label{AS.5} For each $i\in \ii{n}$ and $k\geq 0$, there exist a constant
  $\varsigma_i>0$ and  a random variable $v_{i,k}$ such that $v_{i,k} \geq
  \varsigma_i$ and $w_{i,k} = (v_{i,k})^{\mu}$ for some $\mu \in (0,1)$. In addition,
  \beqn{errvik}
  |\Econd{v_{i,k}}{k} - v_{i,k}| \leq \kappa_v  ( \Econd{g_{i,k}^2}{k} + g_{i,k}^2  )
  \eeqn
  for some $\kappa_v >0$ and all $k\geq0$.
\end{assumption}
\begin{assumption}\label{AS.6}
 For every realization of the algorithm, we have 
    that $ g_{i,k}^2 \leq v_{i,k}$ for all $i\in\ii{n}$ and all $k\geq 0$.
\end{assumption}

\noindent
We immediately note that Assumption~\ref{AS.5} implies that
\beqn{vbounded}
 v_{i,k} \geq \min_{i\in\ii{n}}\varsigma_i \eqdef \varsigma_{\min}
\eeqn
and Assumption~\ref{AS.6} ensures that
\beqn{Eg2vsEv}
\Econd{g_{i,k}^2}{k} \leq \Econd{v_{i,k}}{k}.
\eeqn

\noindent
The first step in our analysis is to derive a parametric bound on the decrease in the
exact linear model of $F$ caused by the step $s_k$, using a technique
inspired by \cite{WardWuBott19} and \cite{defossez2022a}.

\begin{lemma}\label{upperbounddecrease}
 Let $s_j^L$ be the step produced at the $j$-th iteration by the \al{ASGRAD}
  algorithm. Suppose also that Assumptions~\ref{AS.4}, \ref{AS.5} and
  \ref{AS.6} hold. Let $G_j$ be the true gradient of $F$ at
  $x_j$. Then, for all $i \in \ii{n}$, 
  \beqn{isdescent}
  \Econd{\gamma_j G_{i,j} s_{i,j}^L}{j}
  \!\leq \! -(1-\frac{\mu}{2})\frac{\gamlow G_{i,j}^2} {(\Econd{v_{i,j}}{j})^{\mu}}
       + 2 \kappa_\Delta\Econd{\frac{g_{i,j}^2}{w_{i,j}^2}}{j},
  \eeqn
  where 
  \beqn{kappaDelta-def}
  \kappa_\Delta
  \eqdef \frac{\mu \kappa_v^2}{\gamlow} \left[ \kappa_g^{2 \mu}
    + \frac{\kappa_g^2}{\varsigma_{\min}^{1-\mu}}
    + \frac{\kappa_g^{4- 2 \mu}}{\varsigma_{\min}^{2-2 \mu}}
    + \kappa_g^{2 -2\mu}\kappa_\mu\right]
\tim{with}
  \kappa_\mu
  \eqdef \frac{1}{\varsigma_{\min}^{1-2\mu}}\indic{\mu < \half}
                              + \kappa_g^{4\mu-2}\indic{\mu \geq \half},
  \eeqn
where $\indic{\calE}$ stands for the indicator function of the event $\calE$.
\end{lemma}

\vspace*{-5mm}
\proof{See Appendix.}

\noindent
This lemma essentially implies that $s^L$ provides a descent direction on the
true $F$ as long as the square of the true gradient's norm remains large
compared with the stepsizes. We also need another result, partly inspired by
\cite{WardWuBott19,defossez2022a}, whose utility will be to bound
the last term on the right-hand side of \req{isdescent}.

\begin{lemma}\label{gen:series}
Let $\{a_k\}_{k\ge 0}$ be a non-negative sequence,
$\alpha > 0$ and define, for each $k \geq 0$,
$b_k = \sum_{j=0}^k a_j$.  Then if $\alpha \neq 1 $,
\beqn{allalpha series-bound}
\sum_{j=0}^k  \frac{a_j}{(\varsigma+b_j)^{\alpha}}
\le \frac{1}{(1-\alpha)} ( (\varsigma + b_k)^{1-  \alpha} - \varsigma^{1-  \alpha} ).
\eeqn
Otherwise (i.e.\ if $\alpha  = 1$) (see Lemma~5.2 in \cite{defossez2022a}),
\beqn{alphasup1series-bound}
\sum_{j=0}^k  \frac{a_j}{\varsigma+b_j}
\le  \log\left(\frac{\varsigma + b_k}{\varsigma} \right).
\eeqn
\end{lemma}
\vspace*{-4mm}
\proof{See Appendix. Note that \req{alphasup1series-bound} is the limit of \req{allalpha series-bound} when $\alpha$
tends to one.}

\noindent
Using both Lemmas~\ref{upperbounddecrease} and \ref{gen:series}, we are now in
position to deduce a first result on the global convergence rate of a
class of \al{ASGRAD} algorithms using specific ``Adagrad-like'' scaling factors
satisfying Assumptions~\ref{AS.5} and \ref{AS.6}.

\begin{theorem}\label{theorem:stochasticallmu}
Suppose that Assumptions~\ref{AS.1}--\ref{AS.4} hold and that the \al{ASGRAD} algorithm is applied to
problem \req{problem} where, for all $k\geq 0$ and all $i\in\ii{n}$,
\beqn{stow-adapallmu}
w_{i,k} = \left(\varsigma + \sum_{\ell=0}^k g_{i,\ell}^2\right)^\mu,
\eeqn
where $\varsigma \in (0,\kappa_g]$ and $\mu\in (0,1)$.
Then the following bounds hold for $\kappa_\Delta$ given in
\req{kappaDelta-def} and
\beqn{kapbox-def}
\kappa_\Box \eqdef \bigfrac{\kappa_g^{2\mu}(4\kappa_\Delta+L)}{(1-\bigfrac{\mu}{2})\gamlow
  }.
\eeqn
\vspace*{-4mm}
\begin{itemize}
\item[(i) ] If $\mu \in (0,\frac{1}{2})$, then
\begin{align}\label{stogradboundmuinfhalf}
\E{\average_{j\in\iiz{k}}\|G_j\|^2}
&\le \frac{2 \kappa_g^{2\mu}}{ (1-\frac{\mu}{2})\gamlow(k+1)^{1-\mu}}\Big[F(x_0)-\Flow\Big] \nonumber \\
& \hspace*{20mm}+ \frac{n\kappa_\Box}{ 1-2\mu  }
 \frac{(\varsigma+\kappa_g^2(k+1))^{1-2\mu}-\varsigma^{1-2\mu}}{(k+1)^{1-\mu}}.
\end{align}
\item[(ii) ]	
If $\mu = \frac{1}{2}$, then
\begin{align}\label{stogradboundmuequalphalf}
\E{\average_{j\in\iiz{k}}\|G_j\|^2}
& \leq \frac{8\kappa_g}{3\gamlow\sqrt{(k+1)}}\Big[F(x_0)-\Flow\Big]
+ n\kappa_\Box
  \frac{\log\left(  1 + (k+1)\bigfrac{\kappa_g^2}{\varsigma}\right)}{\sqrt{(k+1)}}.
\end{align}
\item[(iii) ]	
If $\mu \in (\frac{1}{2},1)$, then
\begin{align}\label{stogradboundmusupphalf}
\E{\average_{j\in\iiz{k}}\|G_j\|^2}
&\le \bigfrac{2\kappa_g^{2\mu}}{ (1-\bigfrac{\mu}{2})\gamlow(k+1)^{1-\mu}}\Big[F(x_0)-\Flow\Big]  \nonumber\\
& \hspace*{20mm}+ \bigfrac{n\kappa_\Box}{2\mu-1}
   \bigfrac{\varsigma^{1-2\mu}\!\!-\!(\varsigma\!+\!\kappa_g^2(k+1))^{1-2\mu}}{(k+1)^{1-\mu}}.
\vspace*{-2mm}
\end{align}
\end{itemize}
\end{theorem}

\bpr
It is clear from \req{stow-adapallmu} that $w_{i,k}\geq \varsigma^\mu$. Moreover,
if we define $v_{i,k} \eqdef \varsigma +  \sum_{\ell=0}^k g_{i,\ell}^2 $, then
we have that $w_{i,k} = v_{i,k}^\mu$, $v_{i,k} \geq g_{i,k}^2$ and 
\[ 
|{\Econd{v_{i,k}}{k}} - v_{i,k}|
= |{\Econd{g_{i,k}^2}{k}} - g_{i,k}^2 |
\leq  \Econd{g_{i,k}^2}{k}  +   g_{i,k}^2 .
\]
Thus the proposed  scaling factors verify Assumptions~\ref{AS.5} and
\ref{AS.6} with $\kappa_v=1$. 

\noindent
Using \req{lipsch-bound}, we derive that  
\[
F(x_{j+1})
\leq F(x_j) + \gamma_j G_j^T s^L_{j} + \frac{L}{2}  \gamma_j^2 \| s_j^L\|^2
\leq F(x_j) + \gamma_j G_j^T s^L_j + \frac{L}{2}  \| s_{j}^L\|^2.
\]

\noindent
Taking the conditional expectation, using Lemma~\ref{upperbounddecrease},
the fact that $v_{i,j} \leq (k+2) \kappa_g^2$ (because we assumed that
$\varsigma \leq \kappa_g$), \req{sL-def-a}, we deduce
that, for $j \in \iiz{k}$,  
\begin{align*}
\Econd{F(x_{j+1})}{j} &\leq {F(x_j)} + 
\sum_{i=1}^n \Econd{\gamma_j G_{i,j} s^{L}_{i,j}}   {j} + \frac{L}{2} \Econd{ \| s_{j}^L\|^2}{j}, \\
&\leq F(x_j) - \sum_{i=1}^n (1-\frac{\mu}{2}) \gamlow
\frac{G_{i,j}^2}{(\Econd{v_{i,j}}{j})^\mu}+2\kappa_\Delta  \Econd{\frac{g_{i,j}^2}{w_{i,j}^2}}{j} +
\frac{L}{2}  \Econd{ \| s_{j}^L\|^2}{j},\\
&\leq F(x_j) - (1-\frac{\mu}{2}) \gamlow \frac{\|G_j\|^2}{\kappa_g^{2\mu}(k+2)^\mu}
   + \left(\frac{L}{2} + 2 \kappa_\Delta\right)\Econd{ \| s_{j}^L\|^2}{j}.
\end{align*}
We may now take the full expectation and sum the previous inequality
for $j \in \iiz{k}$ to derive that
\begin{align}\label{telescopeeq}
\E{F(x_{{k+1}})}
&\leq F(x_0) - (1-\frac{\mu}{2}) \frac{\gamlow}{\kappa_g^{2\mu}(k+2)^\mu} \sum_{j=0}^{k} \E{\| G_j\|^2}
  + \left(\frac{L}{2} + 2 \kappa_\Delta\right)\sum_{j=0}^{k} \Econd{\| s_{j}^L\|^2}{} \nonumber \\ 
&\leq F(x_0) - (1-\frac{\mu}{2}) \frac{\gamlow}{\kappa_g^{2\mu}(k+2)^\mu} \sum_{j=0}^{k} \E{\| G_j\|^2}
  + \left(\frac{L}{2} + 2 \kappa_\Delta\right) \sum_{i=1}^n\sum_{j=0}^{k} \E{(s_{i,j}^L)^2}. \nonumber \\ 
\end{align} 
Using now Lemma~\ref{gen:series} with $\alpha=2\mu$ for each $s^L_{i,j}$,
\req{sL-def-a}, \req{stow-adapallmu} and Assumption~\ref{AS.4},
we derive that, for $ \mu \in (0, \frac{1}{2})$,
\begin{align*}
\sum_{j=0}^{k}  (s_{i,j}^L)^2
& = \sum_{j=0}^{k}  \frac{g_{i,j}^2}{(\varsigma + \sum_{j=0}^k g_{i,j}^2)^{2\mu}}\\
&\leq \frac{1}{1-2\mu} \left[\left(\varsigma +\sum_{j=0}^k g_{i,j}^2 \right)^{1-2\mu}-\varsigma^{1-2\mu}\right]\\
&\leq \frac{1}{1-2\mu}\left[\Big(\varsigma +(k+1) \kappa_g^2\Big)^{1-2\mu}-\varsigma^{1-2\mu}\right].
\end{align*}
Plugging this inequality in \req{telescopeeq} and using Assumption~\ref{AS.3}, we obtain that 
\begin{align*}
\Flow \leq \E{F(x_{k+1})}
&\leq F(x_0)-(1-\frac{\mu}{2})\frac{\gamlow}{\kappa_g^{2\mu}(k+2)^\mu}\sum_{j=0}^{k} \E{\| G_j\|^2} \\
& \hspace*{15mm} + \frac{n}{1-2\mu}\left(\frac{L}{2} + 2 \kappa_\Delta\right)
   \left[(\varsigma + (k+1) \kappa_g^2)^{1-2\mu}-\varsigma^{1-2\mu}\right]
\end{align*}
and thus, since $(k+2)^\mu \leq 2 (k+1)^\mu$, that 
\begin{align}\label{stousefulineq}
(k+1) \E{\average_{j\in\iiz{k}} \| G_j\|^2}
& \leq \sum_{j=0}^{k} \E{\| G_j\|^2}\\
&\leq \bigfrac{ 2 \kappa_g^{2\mu}(F(x_0) - \Flow)}{(1-\frac{\mu}{2})\gamlow  (k+1)^{-\mu} }  \\
  & \hspace*{5mm}+ \frac{n\left[(\varsigma + \kappa_g^2(k+1))^{1-2\mu}-\varsigma^{1-2\mu}\right]}{(1-2\mu)(k+1)^{-\mu}}
  \left(\frac{\kappa_g^{2\mu}\Big(L+4\kappa_\Delta\Big)}{\gamlow (1-\frac{\mu}{2})}\right), \nonumber 
\end{align}
which is \req{stogradboundmuinfhalf}.

\noindent
If $\mu = \frac{1}{2}$, we reuse (\ref{telescopeeq}) and
Lemma~\ref{gen:series} for each $s^L_{i,j}$ with $\alpha = 1$,  and derive
that, in this case,
\[ 
\E{F(x_{k+1})}
\leq F(x_0)
-\frac{3}{4}\frac{\gamlow}{\sqrt{(k+2)}\kappa_g}\sum_{j=0}^{k}\E{\| G_j\|^2}
+ n \left(\frac{L}{2} + 2 \kappa_\Delta\right)  \log\left(1+ (k+1)\bigfrac{\kappa_g^2}{\varsigma}\right).
\]
By a reasoning similar to that leading to \req{stousefulineq} we now obtain that
\begin{align*}
(k+1)\E{\average_{j\in\iiz{k}} \| G_j\|^2}
&\leq \sum_{j=0}^{k} \E{\| G_j\|^2} \\
&\leq \left(\frac{4}{3}\right) \bigfrac{ 2 \kappa_g(F(x_0) - \Flow) \sqrt{(k+1)}}{\gamlow   } \nonumber  \\
&\hspace*{15mm}+ \left(\frac{4n}{3}\right) \frac{\kappa_g}{\gamlow } ( L + 4 \kappa_\Delta) 
 \log\Big( 1+ (k+1)\frac{\kappa_g^2}{\varsigma} \Big) \sqrt{(k+1)}. \nonumber
\end{align*}
Rearranging the terms yields \req{stogradboundmuequalphalf}. 

\noindent
Finally, if $ \mu \in (\frac{1}{2},1)$, we again reuse (\ref{telescopeeq}) and
Lemma~\ref{gen:series} for each $s^L_{i,j}$ with $\alpha = 2 \mu > 1$, and
deduce that
\begin{align*}
\E{F(x_{k+1})}
& \leq F(x_0) - (1-\frac{\mu}{2}) \frac{\gamlow}{(k+2)^\mu
  \kappa_g^{2\mu}} \sum_{j=0}^{k} \E{\| G_j\|^2}\\
& \hspace*{10mm} + \left(\frac{L}{2} + 2 \kappa_\Delta\right) \frac{n}{2\mu-1}
   \left(\varsigma^{1-2\mu}-(\varsigma+\kappa_g^2(k+1))^{1-2\mu} \right).
\end{align*}
Following the same argument as above yields that
\begin{align*}
(k+1) \E{\average_{j\in\iiz{k}} \| G_j\|^2}
&\leq \sum_{j=0}^{k} \E{\| G_j\|^2} \\
&\leq \bigfrac{ 2 \kappa_g^{2\mu}(F(x_0) - \Flow)}{(1-\frac{\mu}{2})\gamlow  (k+1)^{-\mu} } 
 + \bigfrac{n}{2\mu-1}\left(\frac{\kappa_g^{2\mu}\Big(L+4\kappa_\Delta\Big)}{\gamlow(1-\frac{\mu}{2})}\right)\times\\
&\hspace*{10mm}  \frac{\varsigma^{1-2\mu}-(\varsigma+\kappa_g^2(k+1))^{1-2\mu}}{(k+1)^{-\mu}}. \nonumber 
\end{align*}
Rearranging the terms gives \req{stogradboundmusupphalf}.\epr

\noindent
Note that the last fractions in the last terms of
\req{stogradboundmuinfhalf} and \req{stogradboundmusupphalf} have been written
in a form stressing the continuity with \req{stogradboundmuequalphalf}, but could
obviously be bounded above by the simpler
\[
\frac{(\varsigma+\kappa_g^2)^{1-2\mu}}{(k+1)^\mu}
\tim { and }
\bigfrac{\varsigma ^{1-2\mu}}{(k+1)^{1-\mu}}
\]
respectively.
 
Theorem~\ref{theorem:stochasticallmu} suggests a few comments. The first is that
\req{stogradboundmuinfhalf}, \req{stogradboundmuequalphalf}
and \req{stogradboundmusupphalf} guarantee the
convergence of the \al{ASGRAD} algorithm with \req{stow-adapallmu} to first-order critical points,
because their right-hand sides all tend to zero when $k$ tends to infinity.
The rate at which this convergence occurs, however, differs  for the three
cases, depending on the parameter $\mu$.  If constants are lumped
into a generic $\calO(\cdot)$ notation and using that
\[
\frac{1}{(k+1)} \E{\sum_{j=0}^{k} \|G_j\|} \leq \frac{1}{\sqrt{k+1}} \E{\sqrt{\sum_{j=0}^{k} \|G_j\|^2}}
\leq  \sqrt{\E{\average_{j \in \iiz{k}} \|G_j\|^2}}
\]
where we used Cauchy Schwartz and Jensen's inequality, we obtain, that
\[
\E{\average_{j\in\iiz{k}}\|G_j\|}
\le
\left\{\begin{array}{ll}
\!\!\!\!\calO\left(\bigfrac{1}{(k+1)^{\half\mu}}\right)       \!\!\!& \!\!\!(\mu \in (0,\half)),\\*[3ex]
\!\!\!\!\calO\left(\bigfrac{\log(k+1)}{(k+1)^\quarter} \right)\!\!\!& \!\!\!(\mu = \half),\\*[3ex]
\!\!\!\!\calO\left(\bigfrac{1}{(k+1)^{\half(1-\mu)}}\right)    \!\!\!& \!\!\!(\mu \in (\half,1)).
\end{array}\right.
\]
Examining these ``$k$-order'' bounds indicates that the best bound is that
corresponding to $\mu = \half$.  This is nothing but the standard Adagrad algorithm. 

\subsection{Comparison with prior work for $\mu = \half$}

To provide more context for the reader and to better locate our
  work within the vast literature dealing with the theoretical analysis
of Adagrad, we now discuss our result for $\mu = \frac{1}{2}$.
We immediately note that our bound in
$\mathcal{O}\left(\frac{\log(k+1)}{\sqrt{(k+1)}}\right)$ on
$\E{\average_{j\in\iiz{k}}\|G_j\|}$ is not new for Adagrad and has
been obtained under various assumptions, be it under the requirement
of gradient boundedness
\cite{WardWuBott19,defossez2022a,kavis2022,ZhouTangYangCaoGu24} such 
as our case,  in the unbounded case
\cite{Fawetal22,Wangetal23,Fawetal23,AttiaKoren23,LiuNguyenetal23},
under various Lipschitz smoothness assumptions, see
\cite{Fawetal23,LiuPanZhang24,Wangetal23} and broader noise
assumptions
\cite{HongLin24,AttiaKoren23,kavis2022,LiuNguyenetal23}. These works
either focused on a component-wise Adagrad like the one we analyzed, or study a
variant called Adagrad-Norm that set a global scaling $w_k \propto \sqrt{\sum_{i=1}^{k}
  \|g_k\|^2}$ for each dimension. The final result
was also presented in high-probability
\cite{LiuNguyenetal23,AttiaKoren23,kavis2022} or in expectation just
as our case \cite{WardWuBott19,defossez2022a,Fawetal22,Fawetal23}. For
a compact summary of the last results for the study of Adagrad
variants, see[Table~1]\cite{HongLin24}.  Our contribution is to show
that the choice  $\mu = \frac{1}{2}$ is optimal within a wide range class of adaptive gradient methods \eqref{stow-adapallmu}.  

We also note that
$\mathcal{O}\left(\frac{\log(k+1)}{\sqrt{(k+1)}}\right)$ differs by a
logarithmic term from the  convergence rate of well-tuned stochastic
first order methods proven in \cite{GhadimiLan13}, and so improving
our current bound on $\E{\average_{j\in\iiz{k}}\|G_j\|}$ may be
possible under additional assumptions, such as a tighter variance bound
on $g_{i,k}$.   

\numsection{A ``divergent series'' class of
  \tal{ASGRAD} algorithms}\label{section:divseries}

One might then wonder if
a class of \al{ASGRAD} algorithms exists where improved asymptotic convergence rate can be achieved. This
section considers two cases of interest, both depending on some
constants $\mu \in (0,1)$ and $\varsigma > 0$.  The first, which we call
\emph{maxgi}, is defined, for some $\alpha > 1$ and
$i\in \ii{n}$, by
\beqn{maxgi-def}
w_{i,k} = \xi_{i,k}(k+1)^\mu
\tim{where}
\xi_{i,k}  = \left\{\begin{array}{ll}
            \varsigma  & \tim{if} k = -1,\\
            |g_{i,k}|   &\tim{if} k \geq 0 \tim{and} |g_{i,k}| \geq \alpha \xi_{i,k-1}\\
            \xi_{i,k-1} & \tim{if} k \geq 0 \tim{and} |g_{i,k}|  < \alpha \xi_{i,k-1}\\
            \end{array}\right. 
\eeqn
The second, called \emph{avrgi},  uses for $i\in \ii{n}$
\beqn{avrgi-def}
w_{i,k} = \xi_{i,k}(k+1)^{\mu} \tim{ where } \xi_{i,k} = \max\left[\varsigma, \avggik{k+1}{k}\right].
\eeqn
Before delving into the analysis, we briefly mention that only
  $\xi_{i,k}$ defined in \eqref{maxgi-def} is a monotonically
  increasing sequence whereas this is not necessarily the case
    for $\eqref{avrgi-def}$.  In both cases, $\xi_{i,k}$ lies in
    the interval $[\varsigma, \kappa_g]$ if Assumption~\ref{AS.4} holds.

\noindent
We first state a crucial decrease result.

\begin{lemma}\label{avrdecrease-l}
Suppose that Assumptions~\ref{AS.1}--\ref{AS.4} hold and that the
\al{ASGRAD} algorithm is applied to problem \req{problem}  with its
scaling factors being defined, for some $\mu \in (0,1)$ by
\req{maxgi-def} or \req{avrgi-def}. Then
\beqn{avrdecrease}
\Econd{-\gamma_{j} \bigfrac{G_{i,j} g_{i,j}}{w_{i,j}} }{j}
\leq - \kappa_1 \bigfrac{G_{i,j}^2}{(j+1)^\mu} + \kappa_2\frac{\Pcond{\calA_j}{j}}{(j+1)^\mu} +
\kappa_3\Econd{\frac{g_{i,j}^2}{w_{i,j}^2}}{j},
\eeqn
where $\calA_j$ denotes the event $\{|g_{i,j}| \geq \alpha
\xi_{i,j-1}\}$ and
\begin{enumerate}
\item if \req{maxgi-def} is used,
\beqn{maxgi-kappas}
\kappa_1 =\frac{\gamlow}{\varsigma},
\ms
\kappa_2 = \frac{4\kappa_g^3}{\varsigma^2}
\tim{ and }
\kappa_3 = 0,
\eeqn
\item if \req{avrgi-def} is used,
\beqn{avrgi-kappas}
\kappa_1 = \frac{\gamlow}{2\varsigma},
\ms
\kappa_2 = 0
\tim{ and }
\kappa_3 = \frac{2\kappa_g^2}{\varsigma\gamlow}.
\eeqn
\end{enumerate}
\end{lemma}

\proof{The proof is somewhat technical and given in Appendix.}

Observe that, although  it is well defined for both cases, the event
$\calA_j$ is only relevant for \emph{maxgi} because $\kappa_2=0$ in
\req{avrgi-kappas}.
We also give an important property for the \emph{maxgi} case.

\begin{lemma}\label{sumPA-l}
Suppose that Assumptions~\ref{AS.1} and \ref{AS.4} hold and that the
\al{ASGRAD} algorithm is applied to problem \req{problem}  with its
scaling factors being defined by \req{maxgi-def} for some $\mu \in
(0,1)$.
Then, for all $k\geq 0$,
\beqn{sumPA}
\E{\sum_{j=0}^k \Pcond{\calA_j}{j}}
\leq\left\lfloor\frac{\log(\kappa_g)-\log(\varsigma)}{\log(\alpha)}\right\rfloor
\eqdef \tau_{\max}.
\eeqn
\end{lemma}

\proof{
For $k\geq 0$, let $\tau_k$ be the number of occurrences of $\calA_j$ for
$j\in\iiz{k}$. We must have that, for all $j$,
\[
\kappa_g \geq \xi_{i,j} \geq \xi_{i,-1}\alpha^{\tau_k} = \varsigma\alpha^{\tau_k}
\]
and hence, for all $k$,
\beqn{betakbound}
\sum_{j=0}^k \indic{\calA_j} = \tau_k \leq \left\lfloor\frac{\log(\kappa_g)-\log(\varsigma)}{\log(\alpha)}
\right\rfloor.
\eeqn
Using this bound and the law of total expectation, we thus obtain that for all $k$,
\[
\E{\displaystyle\sum_{j=0}^k\Pcond{\calA_j}{j}}
= \displaystyle\sum_{j=0}^k\E{\Econd{\indic{\calA_j}}{j}}
= \sum_{j=0}^k \E{\indic{\calA_j}}
= \E{\sum_{j=0}^k\indic{\calA_j}}
\leq \left\lfloor\frac{\log(\kappa_g)-\log(\varsigma)}{\log(\alpha)}\right\rfloor.
\]
\epr
}

\begin{theorem}\label{theorem:increase}
Suppose that Assumptions~\ref{AS.1}--\ref{AS.4} hold and that the
\al{ASGRAD} algorithm is applied to problem \req{problem}  with its
scaling factors being defined, for some $\mu \in (0,1)$ by
\req{maxgi-def} or \req{avrgi-def}.
Then, for $\mu \neq \half$,
\beqn{Egbdnu}
\begin{array}{lcl}
\E{\average_{j\in\iiz{k}}\|G_j\|^2}
&\leq &\bigfrac{F(x_0)+n\kappa_2\tau_{\max}- \Flow}{\kappa_1 (k+1)^{1-\mu}}\\*[2ex]
&&+\bigfrac{n\kappa_g^2}{\kappa_1\varsigma^2(1-2\mu)}\left[\kappa_3+\bigfrac{L}{2}\right]\left[\bigfrac{1}{(k+1)^\mu}-\bigfrac{2\mu}{(k+1)^{1-\mu}}\right],
\end{array}
\eeqn
\vspace*{-2mm}
while, if $\mu = \half$,
\begin{align}\label{Egbdhalf}
\E{\average_{j\in\iiz{k}}\|G_j\|^2}
&\leq \frac{F(x_0)+n\kappa_2\tau_{\max}-\Flow}{\kappa_1 \sqrt{k+1}}
   +\bigfrac{n\kappa_g^2}{\kappa_1\varsigma^2}\left[\kappa_3+\bigfrac{L}{2}\right]\frac{1+\log(k+1)}{\sqrt{k+1}},
\end{align}
where $\kappa_1$, $\kappa_2$, $\kappa_3$ and $\tau_{\max}$ are defined
in Lemmas~\ref{avrdecrease-l} and \ref{sumPA-l}.
\end{theorem}

\proof{
By using \req{lipsch-bound}, the inequality $\gamma_j\leq 1$,
and \req{sL-def-a}, we derive that 
\beqn{aa}
F(x_{j+1})
\leq F(x_j) + \gamma_{j} G_j^T s^L_{j} + \frac{L}{2}  \gamma_{j}^2 \| s_j^L\|^2
\leq F(x_j) - \gamma_{j} \frac{G_j^Tg_{i,j}}{w_{i,j}} + \frac{L}{2} \sum_{i=1}^n\frac{g_{i,j}^2}{w_{i,j}^2}.
\eeqn
Using Lemma~\ref{avrdecrease-l}, taking the conditional expectation of
\req{aa} and using Assumption~\ref{AS.4} to bound $g_{i,j}^2$ and that $w_{i,k} \geq \varsigma (k+1)^\mu$ for both choices \eqref{maxgi-def} and \eqref{avrgi-def}, we obtain that
\begin{align}
\Econd{F(x_{j+1})}{j}
&\leq F(x_j)+\sum_{i=1}^n\Econd{\gamma_jG_{i,j}\frac{g_{i,j}}{w_{i,j}}}{j}+\frac{L}{2} \Econd{\frac{g_{i,j}^2}{w_{i,j}^2}}{j}, \nonumber  \\
&\leq F(x_j)-\sum_{i=1}^n\left[\kappa_1\frac{G_{i,j}^2}{(j+1)^\mu}
  +\kappa_2 \frac{\Pcond{\calA_j}{j}}{(j+1)^\mu}+\kappa_3\Econd{\frac{g_{i,k}^2}{w_{i,j}^2}}{j}
   +\frac{L}{2} \Econd{\frac{g_{i,j}^2}{w_{i,j}^2}}{j}\right] \label{futur}  \\
&\leq F(x_j)-\sum_{i=1}^n\frac{\kappa_1 G_{i,j}^2}{(j+1)^\mu}
   + n\kappa_2\Pcond{\calA_j}{j}
   +\frac{n\kappa_g^2}{\varsigma^2}\left[\kappa_3+\frac{L}{2}\right]\frac{1}{(j+1)^{2\mu}}
   \nonumber
\end{align}
Summing over all iterations from 0 to $k$, taking the full expectation
and using Lemma~\ref{sumPA-l} gives that
\[
\begin{array}{lcl}
\E{F(x_{k+1})}
&\leq &F(x_0)- \kappa_1\bigsum_{j=0}^k\bigsum_{i=1}^n\bigfrac{\E{G_{i,j}^2}}{(j+1)^\mu}
+ n\kappa_2 \E{\bigsum_{j=0}^k\Pcond{\calA_j}{j}}
+  \frac{n\kappa_g^2}{\varsigma^2}\left[\kappa_3+\frac{L}{2}\right]\bigsum_{j=0}^k\bigfrac{1}{(j+1)^{2\mu}}\\
&\leq &F(x_0)+n\kappa_2\tau_{\max} - \kappa_1\bigsum_{j=0}^k\bigsum_{i=1}^n\bigfrac{\E{G_{i,j}^2}}{(j+1)^\mu}
+  \frac{n\kappa_g^2}{\varsigma^2}\left[\kappa_3+\frac{L}{2}\right]\bigsum_{j=0}^k\bigfrac{1}{(j+1)^{2\mu}}\\
\end{array}
\]
If we now define
\[
\phi_\mu(x) \eqdef \left\{\begin{array}{ll}
\bigfrac{(x+1)^{1-2\mu}-1}{1-2\mu} & \tim{if } \mu \neq \half\\
\log(x+1)                        & \tim{otherwise,}
\end{array}\right.
\]
we may bound the last inequality, using a simple sum-integral comparison and
Assumption~\ref{AS.3} to obtain that
\[
\sum_{j=0}^k\sum_{i=1}^n\E{G_{i,j}^2}
\leq \frac{(k+1)^\mu(F(x_0)+n\kappa_2\tau_{\max}- \Flow)}{\kappa_1}
   +  \frac{n\kappa_g^2}{\varsigma^2}\left[\kappa_3+\frac{L}{2}\right]\frac{(k+1)^\mu}{\kappa_1}(1+\phi_\mu(k))
\]
and thus that
\[
\E{\average_{j\in\iiz{k}}\|G_j\|^2}
\leq \frac{F(x_0)+n\kappa_2\tau_{\max}- \Flow}{\kappa_1(k+1)^{1-\mu}}
+ \frac{n\kappa_g^2}{\kappa_1\varsigma^2}\left[\kappa_3+\frac{L}{2}\right]\frac{1+\phi_\mu(k)}{(k+1)^{1-\mu}}.
\]
This gives \req{Egbdhalf} when $\mu = \half$.  Otherwise, \req{Egbdnu}
follows from the fact that
\[
1+\phi_\mu(k)
= \frac{1}{1-2\mu}\left[\frac{1}{(k+1)^{2\mu -1}} - 2\mu\right].
\]
\epr
}

\noindent
The choices \req{maxgi-def} and \req{avrgi-def} are of course
reminiscent, in a smooth but stochastic and nonconvex setting, of the
``divergent stepsize'' subgradient method for non-smooth convex
optimization (see \cite{Beck17} and the many references therein), for
which a $\calO(1/\sqrt{k})$ global rate of convergence is known
(Theorems~8.13 and 8.30 in this last reference). 

The bounds given by Theorem~\ref{theorem:increase} are qualitatively similar
to those of Theorem~\ref{theorem:stochasticallmu}, but they may be improved
if we strengthen our assumptions, and impose an additional conditional variance condition on the gradient estimator.

\begin{theorem}\label{theorem:stochasticmemory}
Suppose that Assumptions~\ref{AS.1}--\ref{AS.4} hold and that an
\al{ASGRAD} algorithm is applied to problem \req{problem} with its
scaling factors being defined \req{maxgi-def} and \req{avrgi-def}.
Suppose also that, for all $i \in \ii{n}$ and all $k\geq 0$ 
\beqn{variance-cond}
\Vcond{g_{i,k}}{k}
= \Econd{g_{i,k}^2 - G_{i,k}^2 }{k}
\leq \kap{var} G_{i,k}^2
\eeqn
holds for some $\kap{var} \geq 0$.
Then, for any $\theta \in (0, \kappa_1)$, 
\beqn{stongkbound-ratememort}
\E{\average_{j\in\iibe{j_\theta+1}{k}} \|G_j\|^2\!}
\!\leq\! \kappa_\#(\theta) \frac{(k+1)^{\mu}}{k-j_\theta}
\!\leq\!\frac{\kappa_\#(\theta)(j_\theta\!+\!2)}{(k+1)^{1-\mu}},
\eeqn
where
\beqn{stokappastar-defmemoire}
\kappa_\#(\theta) \eqdef
\frac{1}{\theta} \left(
F(x_0) +n \kappa_2 \tau_{\max} - \Flow +
\frac{n\,2^\mu \kappa_g^4 }{\varsigma^4}\left[\kappa_3+\frac{L}{2}\right] \,(1+\kap{var})j_\theta 
\right),
\eeqn
and
\beqn{defjthetasto}
j_\theta \eqdef
\left\lfloor\left(\left[\kappa_3+\frac{L}{2}\right]\frac{\kappa_g^2 2^\mu}{\varsigma^4(\kappa_1-\theta)}\right)^\sfrac{1}{\mu}\right\rfloor +1.
\eeqn
\end{theorem}

\proof{
To simplify notation, set, for the course of this proof,
$w_{i,-1} = \varsigma$, $i \in \ii{n}$, $\frac{0}{0}=1$.
As in the proof of Theorem~\ref{theorem:increase}, we derive (see \req{futur})
that
\[
\begin{array}{lcl}
\Econd{F(x_{j+1})}{j}
&\leq & F(x_j)-\bigsum_{i=1}^n\left[\kappa_1\bigfrac{G_{i,j}^2}{(j+1)^\mu}
  +\kappa_2 \bigfrac{\Pcond{\calA_j}{j}}{(j+1)^\mu}
  +\left[\kappa_3 +\bigfrac{L}{2}\right] \Econd{\bigfrac{g_{i,j}^2}{w_{i,j}^2}}{j}\right]\\
&\leq & F(x_j)-\bigsum_{i=1}^n\left[\bigfrac{\kappa_1 G_{i,j}^2}{(j+1)^\mu}
  +\kappa_2 \bigfrac{\Pcond{\calA_j}{j}}{(j+1)^\mu}
      +\left[\kappa_3+\bigfrac{L}{2}\right]\Econd{\bigfrac{g_{i,j}^2}{w_{i,j-1}^2}
      \left(\bigfrac{w_{i,j-1}}{w_{i,j}}\right)^2}{j}\right] \nonumber \\
&\leq & F(x_j)-\bigsum_{i=1}^n\left[\bigfrac{\kappa_1 G_{i,j}^2}{(j+1)^\mu}
  +\kappa_2 \bigfrac{\Pcond{\calA_j}{j}}{(j+1)^\mu}
      +\left[\kappa_3+\bigfrac{L}{2}\right]\bigfrac{\kappa_g^2}{\varsigma^2}\Econd{\frac{g_{i,j}^2}{w_{i,j-1}^2} }{j}\right]\nonumber\\
&\leq & F(x_j)-\bigsum_{i=1}^n\left[\bigfrac{\kappa_1 G_{i,j}^2}{(j+1)^\mu}
  +\kappa_2 \Pcond{\calA_j}{j}
      +\left[\kappa_3+\bigfrac{L}{2}\right]\bigfrac{\kappa_g^2(1+\kap{var})}{\varsigma^2w_{i,j-1}^2} {G_{i,j}^2}\right], \nonumber
\end{array}
\]
where we have used the fact that
$\left(\frac{w_{i,j-1}}{w_{i,j}}\right)^2\leq\frac{\kappa_g^2}{\varsigma^2}$ (because of
~\req{maxgi-def} and \req{avrgi-def}),
the measurability of $w_{i,j-1}$  with respect to the past and \req{variance-cond} to
deduce the last inequality.
Using now the bound
$\frac{(j+1)^\mu}{w_{i,j-1}} \leq \frac{2^\mu}{\varsigma}$
and summing over the iterations for $j \in \iiz{k}$ then yields that
\beqn{sumeqstomemory}
\sum_{j=0}^k \Econd{F(x_{j+1})}{j}
\leq \sum_{j=0}^k F(x_j)  +n \kappa_2 \sum_{j=0}^k\Pcond{\calA_j}{j}+ \sum_{j=0}^k \sum_{i=1}^n \frac{G_{i,j}^2}{(j+1)^\mu}
\left(-\kappa_1+\frac{\widehat{\kappa}}{{w_{i,j-1}}} \right)
\eeqn
with $\widehat{\kappa}
= \left[\kappa_3+\frac{L}{2}\right]\frac{\kappa_g^2 2^\mu}{\varsigma^3} (1+\kap{var})$.
Note now that the definition of $j_{\theta}$ in \eqref{defjthetasto} and the fact
that $w_{i,j-1} \geq \varsigma \,j^\mu$ together imply that
\beqn{memoryboundjtheta}
\left(-\kappa_1  +  \frac{\widehat{\kappa}}{w_{i,j-1}}  \right) \leq -\theta,
\eeqn
for $j \geq j_{\theta}$. Hence, from \req{sumeqstomemory},
\beqn{bb}
\begin{array}{ll}
\sum_{j=0}^k \Econd{F(x_{j+1})}{j}
\leq & \bigsum_{j=0}^k F(x_j) +n\kappa_2 \bigsum_{j=0}^k\Pcond{\calA_j}{j}
-\theta \bigsum_{j=j_{\theta}}^k \bigsum_{i=1}^n\bigfrac{G_{i,j}^2}{(j+1)^\mu}\\
 & + \bigsum_{j=0}^{j_\theta-1}\bigsum_{i=1}^n\bigfrac{G_{i,j}^2}{(j+1)^\mu}
\left(-\kappa_1+\bigfrac{\widehat{\kappa}}{w_{i,j-1}} \right),
\end{array}
\eeqn
and the last term of this inequality is bounded by
\begin{align}\label{boundmemory}
\sum_{j=0}^{j_{\theta}-1} \sum_{i=1}^n \frac{G_{i,j}^2}{ (j+1)^\mu}
    \left(-\kappa_1+\frac{\widehat{\kappa}}{{w_{i,j}}} \right)
\leq \sum_{j=0}^{j_{\theta}-1} \sum_{i=1}^n \widehat{\kappa} \frac{G_{i,j}^2}{ \varsigma}
\leq \frac{n \kappa_g^2\widehat{\kappa}}{\varsigma} j_\theta,
\end{align}
where we used the facts that $\| G\|_{\infty} \leq \kappa_g$ (because of
\req{Gbounded}), $w_{i,j}\geq \varsigma$ (because of
\req{maxgi-def} and \req{avrgi-def}).
Injecting \req{boundmemory} in \req{bb}, we deduce that
\[
\theta \sum_{j=j_{\theta}}^k \sum_{i=1}^n \frac{G^2_{i,j}}{(j+1)^\mu}
\leq \sum_{j=0}^k F(x_j) +\kappa_2 \sum_{j=0}^k\Pcond{\calA_j}{j}- \sum_{j=0}^k \Econd{F(x_{j+1})}{j}
     + \frac{n \kappa_g^2  \widehat{\kappa}}{\varsigma} j_\theta.
\]
Taking the full expectation and using Lemma~\ref{sumPA-l} whenever
$\kappa_2 >0$ (i.e.\ for \emph{maxgi}) then gives that
\beqn{lasteqmemory}
(k-j_\theta)\E{\average_{j\in\iibe{j_\theta+1}{k}} \|G_j \|^2 }
\!\leq\! \E{\bigsum_{j=j_\theta}^k \bigsum_{i=1}^n G_{i,j}^2}
\!\leq\! \bigfrac{ (k+1)^{\mu}}{\theta} \left[ F(x_0) +n\kappa_2\tau_{\max}- \Flow
    +   \bigfrac{n \kappa_g^2 \widehat{\kappa}}{\varsigma} j_\theta \right].
\eeqn
which gives the desired result.
\epr
}

\noindent
The (asymptotic) $k$-order of convergence  of
$\E{\average_{j\in\iibe{j_\theta+1}{k}} \|G_j\|}$ implied by
\req{stongkbound-ratememort} is  
therefore
\[
\calO\left(\frac{1}{(k+1)^{\half(1-\mu)}}\right)
\]
where $j_\theta$ is given by \eqref{defjthetasto}.

\numsection{Numerical illustration}\label{section:numerics}

We now provide some numerical illustrations of the algorithmic variants
discussed in the previous sections. We trained
a simple convolutional network  of \cite{GitmGins17} (denoted in
the paper as cifar10-nv) and a small resnet18 model \cite{HeZhanRenSun15} on the CIFAR-10 image
classification dataset\footnote{https://www.cs.toronto.edu/~kriz/cifar.html}.
For these experiments, we used haiku \cite{HennCaiNormBabu20} and optax
\cite{Hessetal20}, two JAX \cite{Bradetal18} based libraries, on a workstation
with four GTX 1080TI. We now compare the numerical performance of
\eqref{stow-adapallmu} for various $\mu$ values in $(0.1,0.5,0.9)$
and of the two scaling factor defined by \req{maxgi-def} and \req{avrgi-def}
with $\mu = 0.1$ $\alpha = 1.1$ and $\varsigma = 0.01$. 
For the experiments, we have chosen two different learning-rate
strategies.
For the cifar10-nv architecture, we chose a fixed\footnote{
Our choice of a fixed learning rate policy is meant to showcase some intrinsic
properties of each scaling factor option.
}
learning rate policy with $\gamma_k = \gamma = 5.\{10^{-4},10^{-5}\}$
for all $k\geq0$. For the resnet18 numerical
  test, we choose a linearly declining learning rate from $\gamma_{\max} = 5. 10^{-2}$ to $\gamma_{\rm low} = 5.10^{-4}$. Note that this choice is covered by our proposed \tal{ASGRAD} framework and is considered as a good choice of schedule for the learning rate as it is built in Optax.
We used the same random initialization for all scaling choices and followed the
data-augmentation procedure of \cite{GitmGins17}, both for training and
testing. 
We trained the models for a
total of $100000$ steps with a batchsize of $128$ using the mean-cross entropy
loss function. We report the training and test accuracies (the
latter on a sample of size $128$ from the test dataset) every $500$ steps.

\begin{figure}[htbp]
  \centerline{\includegraphics[width=0.6\textwidth]{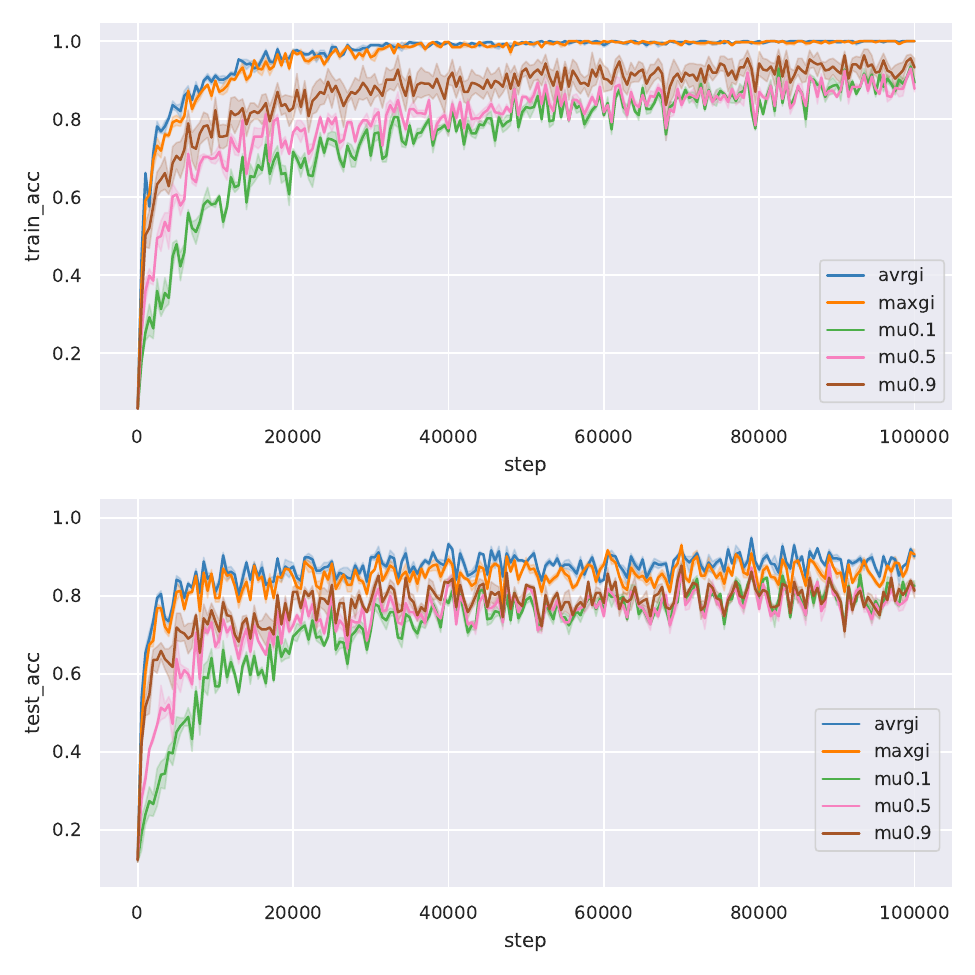}}
\vspace*{-4mm}
\caption{\label{ce-4}Training (top) and test (bottom) accuracies for the
  Adagrad-like ($\mu\in (0.1,0.5,0.9)$), \textit{maxgi} and
  \textit{avrgi} variants with $\gamma=5.10^{-4}$ on the cifar10-nv architecture}
\end{figure}

\begin{figure}[htbp]
  \centerline{\includegraphics[width=0.6\textwidth]{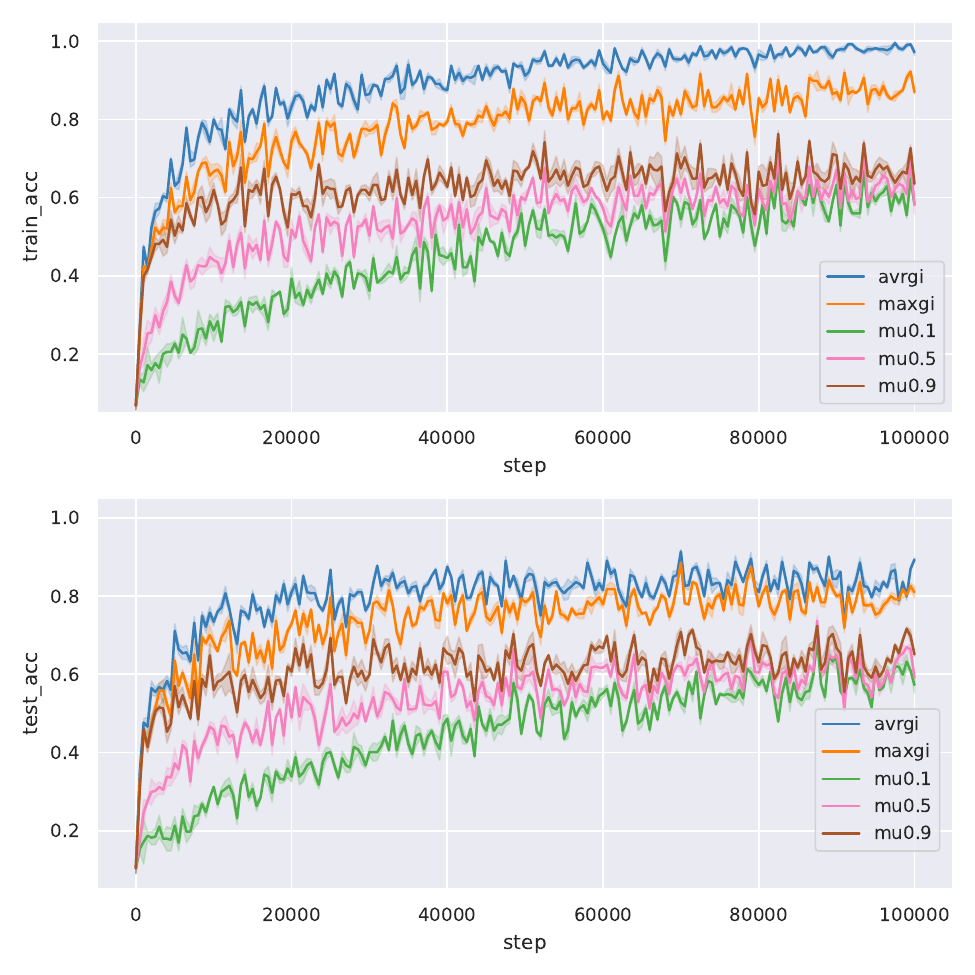}}
\vspace*{-4mm}
\caption{\label{ce-5}Training (top) and test (bottom) accuracy for the
  Adagrad-like ($\mu\in (0.1,0.5,0.9)$), \textit{maxgi} and
  \textit{avrgi} variants with $\gamma=5.10^{-5}$ on the cifar10-nv architecture}
\end{figure}

\begin{figure}[htbp]
  \centerline{\includegraphics[width=0.6\textwidth]{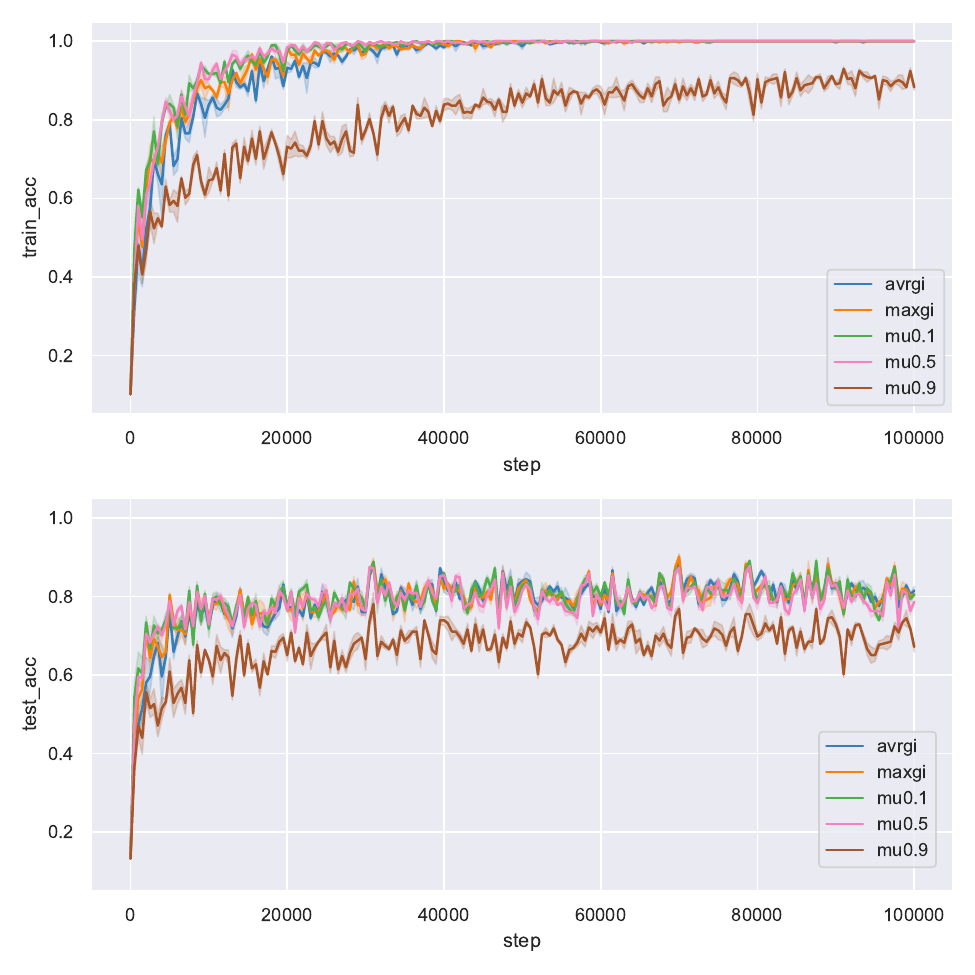}}
\vspace*{-4mm}
\caption{\label{re}Training (top) and test (bottom) accuracies for the
  Adagrad-like ($\mu\in (0.1,0.5,0.9)$), \textit{maxgi} and
  \textit{avrgi} variants with linearly decaying $\gamma$ on the resnet18 architecture}
\end{figure}

The results of these experiments (averaged over three random runs) are
presented in Figures~\ref{ce-4}--\ref{re}.  In each figure, the top panel shows
the evolution (as a function of the number of steps) of the training accuracy,
and the bottom panel that of the test accuracy. The average values are
shown as thick lines and the shaded areas of corresponding colour give the
67\% confidence intervals.


These simple numerical illustrations are obviously not meant to replace significant
numerical testing, but, albeit caution must be exercised not to extrapolate
from limited data, they still suggest a few tentative comments.

\vspace*{-2mm}
\begin{itemize}
\item The relative behaviour of the tested variants  differs
  significantly between the two tested network architectures, even if the
  test accuracy is (as expected) slighly lower for the resnet18 case.
\vspace*{-2mm}
\item For fixed learning rates, the methods \textit{maxgi} and \textit{avrgi}
  of the second \al{ASGRAD} class (introduced in
  Section~\ref{section:divseries})  seem to produce relatively good results on
  our example for fixed learning rate and for the civar-nv architecture, both in training and testing, often outperforming the
  Adagrad-like variants of the first class (of Section~\ref{section:adag}).
\vspace*{-2mm}
\item Among Adagrad-like variants, those with a larger $\mu$ handle
  smaller and fixed
  learning rates  better on these examples, a behaviour admittedly not
  predicted by our theory.
\item The choice of a learning rate schedule for a
  specific architecture has an impact in practice. We see that for the
  resnet18 architecture (Figure~\ref{re}), all the methods behave
  very comparably (except for one variant).
\vspace*{-2mm}
\item The comparison of Figures~\ref{ce-4}, \ref{ce-5} and
  Figure~\ref{re}  unsurprisingly shows that, albeit our
  theory does not depend on the choice of $\gamma_k$, the practical
  convergence behaviour may be affected by this choice (and other factors such
  as the batchsize).
\end{itemize}
\vspace*{-2mm}

\numsection{Conclusions}\label{section:conclusion}

We have introduced a first-order trust-region framework for minimization
methods and derived complexity upper bounds for two classes of interest,
the first containing the standard Adagrad. These bounds give the best
complexity to values of the class parameters corresponding to Adagrad
in the first class. We have also shown these bounds can be improved for both
classes under an additional variance condition, in which case the parameter
choice yielding the best bounds no longer corresponds to Adagrad. This
improvement is asymptotic and implicit for the first class and explicit for
the second. However, our numerical illustrations of the discussed methods on
examples arising from deep-learning applications indicate that methods
of the second class have merits, but also that, at least in our
examples, there remains some distance from the above theory to real
behaviour. This may possibly be because the complexity bounds may not
be sharp, but also, fortunately, because the worst-case happens very
rarely in practice.

{\footnotesize
\paragraph{Acknowledgments}
This work was supported in part by
3IA Artificial and Natural Intelligence Toulouse Institute,
French ”Investing for the Future - PIA3” program under the
Grant agreement ANR-19-PI3A-0004”. The experiments
presented in this paper were carried out using the OSIRIM
platform that is administered by IRIT and supported by
CNRS, the Region Midi-Pyr\'en\'ees, the French Government,
and ERDF (see http://osirim.irit.fr/site/en).


}

\newpage
\appendix
\onecolumn

\appnumsection{A first technical lemma}

\begin{lemma}\label{rupperbound}
Let $\mu \in (0,1]$. Let $x, \,y \in \mathbb{R}^+\setminus\{0\}$.
Then
\beqn{rmajor}
\frac{|x^\mu - y^\mu|}{x^\mu y^\mu} \leq \mu \frac{|x-y|}{x y^\mu} + \mu \frac{|x-y|}{x^\mu y}.
\eeqn 
\end{lemma}

\proof{
Let us first consider the case $x\geq y$. 
Remembering that $u^\mu \leq 1 + \mu (u-1)$ for $u >0$ and taking
$u=\frac{x}{y}$, we successively derive that
\begin{align}\label{xsupy}
     \frac{x^\mu}{y^\mu} &\leq  1 + \mu \left(\frac{x}{y} - 1\right), \nonumber \\
     x^\mu - y^\mu &\leq \mu \left(\frac{x y^\mu}{y} - y^\mu \right)
      =  \mu y^{\mu-1} (x - y ), \nonumber \\
     \frac{x^\mu - y^\mu}{x^\mu y^\mu} &\leq \mu \frac{x-y}{x^\mu y}.
\end{align} 
Hence the inequality \eqref{rmajor} is valid when $x\geq y$. 
For the symmetric case ($y \geq x$), we similarly obtain that
\beqn{ysupx}
\frac{y^\mu - x^\mu}{x^\mu y^\mu} \leq \mu \frac{y-x}{y^\mu x}.
\eeqn
Combining  \eqref{xsupy} and \eqref{ysupx} yields the desired result.
}

\appnumsection{Proof of Lemma~\ref{upperbounddecrease}}

Let us consider an iteration index $j\geq 0$ and a component index
$i\in\ii{n}$. 
We first use the definition of $s^L$ in \req{sL-def-a} and the fact that $w_{i,j} = v_{i,j}^\mu$  (Assumption~\ref{AS.5}) to obtain that
\beqn{eq1}
\Econd{\gamma_j G_{i,j} s^L_{i,j}}{j} = -  \Econd{ \gamma_j  \frac{G_{i,j} g_{i,j}}{v_{i,j}^{\mu}}}{j}
= - \Econd{ \gamma_j \frac{G_{i,j} g_{i,j}}{\Econd{v_{i,j}}{j}^{\mu}}}{j}
  + \Econd{  \gamma_j G_{i,j} g_{i,j} \left( \frac{1}{\Econd{v_{i,j}}{j}^{\mu}}
  - \frac{1}{v_{i,j}^{\mu}}\right)}{j}.
\eeqn
Using that $G_{i,j}$ and $\gamma_{j}$ are
mesurable with respect to the past and  Assumption~\ref{AS.4}, we derive that, 
\beqn{ineq2}
\Econd{-\frac{\gamma_j G_{i,j}g_{i,j}}{\Econd{v_{i,j}}{j}^\mu} }{j}
= - \frac{\gamma_{j}G_{i,j}}{\Econd{v_{i,j}}{j}^\mu} \Econd{g_{i,j}}{j}
= - \frac{ \gamma_{j} G_{i,j}^2}{\Econd{v_{i,j}}{j}^\mu}
\leq -\gamlow \frac{G_{i,j}^2}{\Econd{v_{i,j}}{j}^\mu},
\eeqn 
where we used the measurability of $\Econd{v_{i,j}}{j}^\mu$ with
respect to the past.
Combining \req{eq1} and \req{ineq2} gives that
\beqn{Atermineq}
\Econd{\gamma_j G_{i,j} s^L_{i,j}}{j}
\leq  -\gamlow \frac{G_{i,j}^2}{\Econd{v_{i,j}}{j}^{\mu}} 
      +\Econd{\underbrace{\gamma_j G_{i,j} g_{i,j} \frac{v_{i,j}^{\mu}-\Econd{v_{i,j}}{j}^{\mu}}{v_{i,j}^{\mu}\Econd{v_{i,j}}{j}^{\mu}}}_A}{j}.
\eeqn
We now derive an upper bound on the absolute value of the $A$ term by
successively using Lemma~\req{rupperbound}, Assumption~\ref{AS.5} and the bound $\gamma_j\leq 1$
 to obtain that
\begin{align*}
|A|
&= |\gamma_j G_{i,j} g_{i,j}| \frac{|v_{i,j}^\mu-\Econd{v_{i,j}}{j}^\mu|}{v_{i,j}^\mu \Econd{v_{i,j}}{j}^\mu} 
\leq   \mu |\gamma_j G_{i,j} g_{i,j}| \frac{|v_{i,j} - \Econd{v_{i,j}}{j} |}{v_{i,j}^\mu \Econd{v_{i,j}}{j}}
+\mu |\gamma_j G_{i,j} g_{i,j}| \frac{|v_{i,j} - \Econd{v_{i,j}}{j}|}{v_{i,j}\Econd{v_{i,j}}{j}^\mu}\\
&\leq \mu \underbrace{|G_{i,j} g_{i,j}| \kappa_v \frac{\Econd{g_{i,j}^2}{j}}{v_{i,j}^\mu \Econd{v_{i,j}}{j}}}_B
    + \mu \underbrace{|G_{i,j} g_{i,j}| \kappa_v    \frac{g_{i,j}^2}{v_{i,j}^\mu \Econd{v_{i,j}}{j}}}_C  \\
&\hspace*{2cm} + \mu \underbrace{|G_{i,j} g_{i,j}| \kappa_v\frac{\Econd{g_{i,j}^2}{j}}{v_{i,j}\Econd{v_{i,j}}{j}^\mu}}_D 
    + \mu \underbrace{|G_{i,j} g_{i,j}| \kappa_v \frac{g_{i,j}^2} {v_{i,j}\Econd{v_{i,j}}{j}^\mu}}_E.
\end{align*}
We now use Young's inequality with $p=q=2$, that is
\beqn{realupper}
\forall \lambda > 0, x, \, y \in \mathbb{R}^+, \, xy \leq \frac{\lambda}{2} x^2 + \frac{y^2}{2 \lambda},
\eeqn
to successively handle the four terms in the last bound.

$\bullet$ For the first term $B$, we choose
\[
x = \frac{|G_{i,j}|}{\Econd{v_{i,j}}{j}^\mu},
\ms
\lambda = \frac{\gamlow \Econd{v_{i,j}}{j}^\mu}{4}
\tim{and}
y =  \kappa_v | g_{i,j} | \frac{\Econd{g_{i,j}^2}{j}}{v_{i,j}^\mu \Econd{v_{i,j}}{j}^{1-\mu}}.
\]
Using \eqref{realupper}, Assumptions~\ref{AS.4},~\ref{AS.6} and \req{Eg2vsEv}, we obtain that
\begin{align*}
B
&\leq \gamlow \frac{G_{i,j}^2}{8 \Econd{v_{i,j}}{j}^\mu}
  + 2 \frac{\kappa_v^2}{\gamlow} \frac{g_{i,j}^2}{v_{i,j}^{2\mu}}\frac{\Econd{g_{i,j}^2}{j}^2 }{\Econd{v_{i,j}}{j}^{2-\mu}} ,\\
&\leq \gamlow \frac{G_{i,j}^2}{8 \Econd{v_{i,j}}{j}^\mu}
  + 2 \frac{\kappa_v^2}{\gamlow}\Econd{g_{i,j}^2}{j}^{ \mu} \frac{g_{i,j}^2}{v_{i,j}^{2 \mu}}\\ 
&\leq \gamlow \frac{G_{i,j}^2}{8 \Econd{v_{i,j}}{j}^\mu}
  + 2 \frac{\kappa_v^2}{\gamlow} \kappa_g^{2 \mu} \frac{g_{i,j}^2}{v_{i,j}^{2 \mu}}. 
\end{align*}
Taking now the expectation over $\Econd{.}{j}$ yields that
\beqn{Bterm}
\Econd{B}{j}
\leq \gamlow \frac{G_{i,j}^2}{8 \Econd{v_{i,j}}{j}^\mu}
  + 2 \frac{\kappa_v^2}{\gamlow} \kappa_g^{2\mu} \Econd{\frac{g_{i,j}^2}{w_{i,j}^2}}{j}.
\eeqn

$\bullet$
Now consider the $C$ term. In this case, we choose
\[
x =\frac{|G_{i,j}g_{i,j}|}{\Econd{v_{i,j}}{j}^\mu},
\ms
\lambda = \gamlow \frac{\Econd{v_{i,j}}{j}^\mu}{4\Econd{g_{i,j}^2}{j}}
\tim{and}
y = \kappa_v \frac{g_{i,j}^2}{v_{i,j}^\mu \Econd{v_{i,j}}{j}^{1-\mu}}
\]
to deduce from \req{realupper} that
\begin{align*}
C
&\leq \gamlow \frac{G_{i,j}^2}{8 \Econd{v_{i,j}}{j}^\mu}\frac{g_{i,j}^2}{\Econd{g_{i,j}^2}{j}}
  + 2  \frac{\kappa_v^2}{\gamlow} \frac{g_{i,j}^4} {v_{i,j}^{2 \mu}}\frac{ \Econd{g_{i,j}^2}{j}}{ \Econd{v_{i,j}}{j}^{2-\mu}} \\ 
&\leq  \gamlow \frac{G_{i,j}^2}{8 \Econd{v_{i,j}}{j}^\mu}\frac{g_{i,j}^2}{\Econd{g_{i,j}^2}{j}}
  + 2  \frac{\kappa_v^2}{\gamlow}  \kappa_g^2 \frac{g_{i,j}^2}{v_{i,j}^{2 \mu}} \frac{1}{\Econd{v_{i,j}}{j}^{1-\mu}} \\ 
&\leq \gamlow \frac{G_{i,j}^2}{8 \Econd{v_{i,j}}{j}^\mu}\frac{g_{i,j}^2}{\Econd{g_{i,j}^2}{j}}
  + 2 \frac{\kappa_v^2}{\gamlow} \frac{\kappa_g^2}{\varsigma_{\min}^{1-\mu}} \frac{g_{i,j}^2}{v_{i,j}^{2 \mu} }, 
\end{align*}
where we successively used the facts that $ \Econd{g_{i,j}^2}{j} \leq
\Econd{v_{i,j}}{j} $ (because of \req{Eg2vsEv}), $g_{i,j}^2 \leq \kappa_g^2$ (because of
Assumption~\ref{AS.4}) and $\Econd{v_{i,j}}{j}^{1-\mu} \geq
\varsigma_{\min}^{1-\mu}$ (because of \req{vbounded}). 
Taking the expectation over $\Econd{.}{j}$ then gives that
\beqn{Cterm}
\Econd{C}{j}
\leq \gamlow\frac{G_{i,j}^2}{8 \Econd{v_{i,j}}{j}^\mu}
  +  2 \frac{\kappa_v^2}{\gamlow}\frac{\kappa_g^2}{\varsigma_{\min}^{1-\mu}} \Econd{\frac{g_{i,j}^2}{w_{i,j}^2}}{j}.
\eeqn
\noindent
(Note that we can divide by $\Econd{g_{i,j}^2}{j}$ above, as it suffice to notice that
$\Econd{g_{i,j}^2}{j} = 0$ implies $g_{i,j}^2 = 0$. $C$ would then be equal to
zero and \eqref{Cterm} would still be verified.)

$\bullet$
Let us now handle the $D$ term. Choosing
\[
x = \frac{|G_{i,j}|}{\Econd{v_{i,j}}{j}^\mu},
\ms
\lambda = \gamlow \frac{\Econd{v_{i,j}}{j}^\mu}{4}
\tim{and}
y =  \kappa_v | g_{i,j} | \frac{\Econd{g_{i,j}^2}{j}}{v_{i,j}},
\]
we now deduce from \req{realupper} that
\begin{align*}
D
&\leq \gamlow \frac{G_{i,j}^2}{8 \Econd{v_{i,j}}{j}^\mu}
 + 2 \frac{\kappa_v^2}{\gamlow} \frac{g_{i,j}^2 \Econd{g_{i,j}^2}{j}^2}{\Econd{v_{i,j}}{j}^\mu v_{i,j}^2}, \\ 
&\leq \gamlow \frac{G_{i,j}^2}{8 \Econd{v_{i,j}}{j}^\mu}
 + 2 \frac{\kappa_v^2}{\gamlow} \frac{g_{i,j}^2}{v_{i,j}^{2\mu}}\frac{1}{v_{i,j}^{2-2 \mu}}
   \frac{\Econd{g_{i,j}^2}{j}^2}{\Econd{v_{i,j}}{j}^\mu}\\
&\leq \gamlow \frac{G_{i,j}^2}{8 \Econd{v_{i,j}}{j}^\mu}
 + 2 \frac{\kappa_v^2}{\gamlow} \frac{g_{i,j}^2}{v_{i,j}^{2\mu}} \frac{1}{v_{i,j}^{2-2 \mu}} \Econd{g_{i,j}^2}{j}^{2-\mu} \\
&\leq \gamlow  \frac{G_{i,j}^2}{8 \Econd{v_{i,j}}{j}^\mu}
 + 2 \frac{\kappa_v^2}{\gamlow} \frac{\kappa_g^{4- 2 \mu}}{\varsigma_{\min}^{2-2 \mu}}\frac{g_{i,j}^2}{v_{i,j}^{2\mu}},
\end{align*}
where, as for the $C$ term, we used the facts that
$\Econd{g_{i,j}^2}{j}^\mu \leq \Econd{v_{i,j}}{j}^\mu$, $g_{i,j}^2 \leq
\kappa_g^2$ and $v_{i,j}^{2-2\mu} \geq \varsigma_{\min}^{2-2\mu}$.
Taking the expectation $\Econd{.}{j}$ yields, in this case, that
\beqn{Dterm}
\Econd{D}{j}
\leq \gamlow \frac{G_{i,j}^2}{8 \Econd{v_{i,j}}{j}^\mu}
  + 2 \frac{\kappa_v^2}{\gamlow} \frac{\kappa_g^{4-2\mu}}{\varsigma_{\min}^{2-2 \mu}}
     \Econd{\frac{g_{i,j}^2}{w_{i,j}^2}}{j}.
\eeqn

$\bullet$
Finally consider the $E$ term. Choosing
\[
x = \frac{|G_{i,j}g_{i,j}|}{\Econd{v_{i,j}}{j}^\mu},
\ms
\lambda = \gamlow \frac{\Econd{v_{i,j}}{j}^\mu}{4\Econd{g_{i,j}^2}{j}}
\tim{and}
y =\kappa_v \frac{g_{i,j}^2}{v_{i,j}}
\]
in \req{realupper} then gives that
\vspace*{-2mm}
\begin{align*}
\Econd{E}{j}
&\leq \gamlow \frac{G_{i,j}^2}{8\Econd{v_{i,j}}{j}^\mu}\frac{g_{i,j}^2}{\Econd{g_{i,j}^2}{j}}
   + 2 \frac{\kappa_v^2}{\gamlow}\frac{g_{i,j}^4 \Econd{g_{i,j}^2}{j}}{\Econd{v_{i,j}}{j}^\mu v_{i,j}^{2}} \\ 
&\leq \gamlow \frac{G_{i,j}^2}{8\Econd{v_{i,j}}{j}^\mu}\frac{g_{i,j}^2}{\Econd{g_{i,j}^2}{j}}
   + 2 \frac{\kappa_v^2}{\gamlow} \Econd{g_{i,j}^2}{j}^{1-\mu}
			\frac{g_{i,j}^2}{v_{i,j}^{2\mu}} \left( 
			\frac{1}{v_{i,j}^{1-2\mu}} \frac{g_{i,j}^2}{v_{i,j}} \indic{\mu < \half} 
		    + \frac{|g_{i,j}^{4-4\mu}|}{v_{i,j}^{2-2\mu}} |g_{i,j}^{4\mu-2}|\indic{\mu \geq \half}\right)   \\ 
&\leq \gamlow \frac{G_{i,j}^2}{8\Econd{v_{i,j}}{j}^\mu}\frac{g_{i,j}^2}{\Econd{g_{i,j}^2}{j}}
    + 2  \frac{\kappa_v^2}{\gamlow} \kappa_g^{2-2\mu} \frac{g_{i,j}^2}{v_{i,j}^{2\mu}} \left( 
		\frac{1}{\varsigma_{\min}^{1-2\mu}} \indic{\mu < \half}
		+  \kappa_g^{4\mu - 2} \indic{\mu \geq \half}	\right), 
\end{align*}
where we once more used the facts that $\Econd{g_{i,j}^2}{j}^\mu \leq
\Econd{v_{i,j}}{j}^\mu$ and $|g_{i,j}| \leq \kappa_g$, in turn implying that
\[
g_{i,j}^2 \leq v_{i,j} \tim{and} v_{i,j} \geq \varsigma_{\min} \tim{ if } \mu < \half
\]
and
\[
|g_{i,j}^{4-4\mu} | \leq v_{i,j}^{2-2\mu} \tim{and} |g_{i,j}^{4\mu-2} | \leq \kappa_g^{4\mu-2} \tim{ if } \mu \geq \half.
\]
Taking the expectation $\Econd{.}{j}$, we deduce that 
\beqn{Eterm}
\Econd{E}{j}
\leq \gamlow \frac{G_{i,j}^2}{8\Econd{v_{i,j}}{j}^\mu}
   + \frac{\kappa_v^2}{\gamlow} \kappa_g^{2 -2\mu}\left( 
      \frac{1}{\varsigma_{\min}^{1-2\mu}} \indic{\mu < \half} 
         + \kappa_g^{4\mu-2} \indic{\mu \geq \half}\right) \Econd{\frac{g_{i,j}^2}{w_{i,j}^2}}{j}.
\eeqn
\vspace*{2mm}
\noindent
Summing now  \eqref{Bterm}, \eqref{Cterm}, \eqref{Dterm} and \eqref{Eterm}
and substituting the obtained upper-bound of $A$ in \eqref{Atermineq}, we
finally obtain \req{isdescent} with \req{kappaDelta-def}.

\appnumsection{Proof of Lemma~\ref{gen:series}}

Consider first the case where  $\alpha \neq 1$ and note that
$\frac{1}{(1-\alpha)} x^{1-\alpha}$ is then a non-decreasing and concave
function on $(0,+\infty)$. Setting $b_{-1} = 0$ and using these properties,
we obtain that, for $j\geq 0$,
\begin{align*}
\frac{a_j}{(\varsigma + b_j)^\alpha} &\leq \frac{1}{1-\alpha}\left( (\varsigma
  + b_j)^{1-\alpha}- (\varsigma + b_j - a_j)^{1-\alpha}\right) \\
&\leq \frac{1}{1-\alpha} \left( (\varsigma + b_j)^{1-\alpha} - (\varsigma + b_{j-1})^{1-\alpha}\right).
\end{align*} 
We then obtain \eqref{allalpha series-bound} by summing this inequality for $j\in\iiz{k}$.

Suppose now that $\alpha = 1$, We then use the concavity and non-decreasing
nature of the logarithm to derive that
\[
\frac{a_j}{(\varsigma + b_j)^\alpha}
= \frac{a_j}{(\varsigma + b_j)}
\leq  \log(\varsigma + b_j) - \log(\varsigma + b_j - a_j) 
\leq  \log(\varsigma + b_j) - \log(\varsigma + b_{j-1}).
\]
The inequality \req{alphasup1series-bound} then again follows by
summing for $j\in\iiz{k}$.

\appnumsection{Proof of Lemma~\ref{avrdecrease-l}}

We have that
\begin{align}\label{decreasbound}
\Econd{-\gamma_{j} \frac{G_{i,j} g_{i,j}}{w_{i,j}} }{j}
&\leq -\gamma_{\text{low}} \frac{G_{i,j}^2}{\Econd{w_{i,j}}{j}} +
\Econd{ |\gamma_{j} G_{i,j} g_{i,j} | \frac{ | w_{i,j} - \Econd{w_{i,j}}{j} | }{w_{i,j}\Econd{w_{i,j}}{j}} }{j} \nonumber\\
&=-\gamma_{\text{low}} \frac{G_{i,j}^2}{\Econd{w_{i,j}}{j}} + \Econd{ |\gamma_{j} G_{i,j} g_{i,j}| \frac{ | \xi_{i,j} - \Econd{\xi_{i,j}}{j} | }{\xi_{i,j} (j+1)^\mu\Econd{\xi_{i,j}}{j}} }{j} \nonumber\\
&\leq  -\gamma_{\text{low}} \frac{G_{i,j}^2}{\Econd{w_{i,j}}{j}} + \Econd{A}{j} 
\end{align}	
where
\beqn{A-def}
A
=  \left|\gamma_{j} G_{i,j} g_{i,j} \right|\frac{ |\xi_{i,j}  - \Econd{\xi_{i,j}}{j} | }{ \xi_{i,j}
  (j+1)^\mu \Econd{\xi_{i,j} }{j} }
\eeqn
In order to compute a bound on the right-hand side of
\req{decreasbound}, we consider the \emph{maxgi} and \emph{avrgi}
cases separately.
Consider the  \emph{maxgi} case first. We see that
\beqn{Exi}
  \begin{array}{lcl}
\Econd{|\xi_{i,j}  - \Econd{\xi_{i,j}}{j} |}{j}
& = & \Econd{\left|\indic{\calA_j}\xi_{i,j}-\Econd{\indic{\calA_j}\xi_{i,j}}{j}
              +\indic{\calA_j^C}\xi_{i,j}-\Econd{\indic{\calA_j^C}\xi_{i,j}}{j} \right|}{j}\\
& = & \Econd{\left|\indic{\calA_j}|g_{i,j}|-\Econd{\indic{\calA_j}|g_{i,j}|}{j}
              +\indic{\calA_j^C}\xi_{i,j-1}-\Econd{\indic{\calA_j^C}\xi_{i,j-1}}{j} \right|}{j}\\
& \leq & \Econd{\left|\indic{\calA_j}|g_{i,j}|-\Econd{\indic{\calA_j}|g_{i,j}|}{j}
              +\xi_{i,j-1}\left(\indic{\calA_j^C}-\Econd{\indic{\calA_j^C}}{j}\right)
              \right|}{j}\\
& \leq & \Econd{\indic{\calA_j}|g_{i,j}|+\Econd{\indic{\calA_j}|g_{i,j}|}{j}
              +\xi_{i,j-1}\left|\indic{\calA_j^C}-\Econd{\indic{\calA_j^C}}{j}
              \right|}{j}\\
&\leq & \kappa_g ( \Econd{\indic{\calA_j}+\Econd{\indic{\calA_j}}{j}}{j})
      +\kappa_g\Econd{\left|\indic{\calA_j^C}-\Econd{\indic{\calA_j^C}}{j}
              \right|}{j}\\
&\leq & 2 \kappa_g\Pcond{\calA_j}{j}
      +\kappa_g\Econd{\left|\indic{\calA_j^C}-\Econd{\indic{\calA_j^C}}{j}
              \right|}{j},
\end{array}
  \eeqn
where we used the bound $\xi_{i,j}\leq \kappa_g$ for all $j$
(resulting from Assumption~\ref{AS.4}) to obtain the penultimate inequality.
Now
\[
\Econd{\left|\mathds{1}_{A_j^C} -
  \Econd{\indic{A_j^C}}{j}\right|}{j}
= \Econd{|1-\indic{A_j} - \Econd{1-\indic{A_j}}{j}|}{j} 
= \Econd{|\indic{A_j} - \Econd{\indic{A_j}}{j}|}{j}
\leq 2 \,\mathbb{P}_j(A_j).
\]
We then obtain \req{avrdecrease}-\req{maxgi-kappas} by substituting this last inequality in
\req{Exi} and combing the result with the bound
\[
\frac{ | \gamma_{j} G_{i,j} g_{i,j}|}
     { \xi_{i,j}(j+1)^\mu \Econd{\xi_{i,j}}{j}  }
\leq \frac{\kappa_g^2}{\varsigma^2(j+1)^\mu},
\]
\req{decreasbound} and \req{A-def}.

\noindent
Now consider the \emph{avrgi} case. Analogously to \req{Exi} and using
the identity $\max(a,b) = \frac{1}{2} (a + b + |a-b|)$, we deduce from
\req{avrgi-def} that
\begin{align*}
|\Econd{\xi_{i,j}}{j} - \xi_{i,j}|
&= \frac{1}{2} \left| \Econd{ \varsigma + \avggik{j+1}{j} +  |
  \avggik{j+1}{j} - \varsigma  |}{j}\right.\\
&\hspace*{10mm} \left. -\varsigma - \avggik{j+1}{j} -  | \avggik{j+1}{j} - \varsigma  | \right| \\
&= \frac{1}{2} \left|   \frac{1}{(j+1)} \Econd{|g_{i,j}|}{j} +
\Econd{| \avggik{j+1}{j} - \varsigma  |}{j} -  \frac{1}{(j+1)}|g_{i,j}| - | \avggik{j+1}{j} - \varsigma  | \right| \\
&\leq \frac{1}{2(j+1)} \bigg\lvert \Econd{|g_{i,j}|}{j} - |g_{i,j}| \bigg\rvert \\
&\ms + \frac{1}{2} \left| \frac{1}{(j+1)}\Econd{|g_{i,j}|}{j} + |\avggik{j+1}{j-1} - \varsigma| + \frac{1}{(j+1)}{|g_{i,j}|}- |\avggik{j+1}{j-1} - \varsigma| \right| \\
&\leq \frac{1}{(j+1)} \bigg\lvert \Econd{|g_{i,j}|}{j} + |g_{i,j}| \bigg\rvert,
\end{align*}
where we used that  $|\avggik{j+1}{j-1} - \varsigma|$ is mesurable
with respect to the past. This inequality, the definition \eqref{A-def} and the bounds $\gamma_{j} \leq 1$ and
$(j+1)^{\mu/2} \leq j+1$ then give that
\begin{align}\label{Abound}
A
\leq \underbrace{ \frac{| G_{i,j}| g_{i,j}^2}{(j+1)^{\mu/2+\mu} \Econd{\xi_{i,j}}{j} \xi_{i,j}}}_B + \underbrace{ \frac{| G_{i,j} g_{i,j} \Econd{|g_{i,j}|}{j}|}{(j+1)^{\mu/2+\mu}  \Econd{\xi_{i,j}}{j} \xi_{i,j}}}_C.
\end{align}
We now use Young's inequality with $p=q=2$, that is
\beqn{realuppercase}
\forall \lambda > 0, x, \, y \in \mathbb{R}^+, \, xy \leq \frac{\lambda}{2} x^2 + \frac{y^2}{2 \lambda},
\eeqn
to successively handle the two terms of \req{Abound}.

$\bullet$ For the first term $B$, we choose
\[
x = \frac{|G_{i,j}|}{(j+1)^{\mu/2} \sqrt{\Econd{\xi_{i,j}}{j}}},
\ms
\lambda = \frac{\gamlow }{2}
\tim{and}
y =   \frac{g_{i,j}^2 }{\sqrt{\Econd{\xi_{i,j}}{j}} \xi_{i,j} (j+1)^\mu}
\]
yielding
\[
B \leq \frac{\gamma_{\text{low}} G_{i,j}^2}{4(j+1)^\mu \Econd{\xi_{i,j}}{j}} + \frac{1}{\gamlow} \frac{g_{i,j}^4}{\Econd{\xi_{i,j}}{j}\xi_{i,j}^2(j+1)^{2\mu}} \\
\leq \frac{\gamma_{\text{low}} G_{i,j}^2}{4(j+1)^\mu \Econd{\xi_{i,j}}{j}} + \frac{1}{\gamlow} \frac{g_{i,j}^2 \kappa_g^2}{\varsigma\xi_{i,j}^2(j+1)^{2\mu}},
\]
where we used that $|g_{i,j}| \leq \kappa_g$ and $\xi_{i,j}\geq \varsigma$. Taking now $\Econd{\cdot}{j}$ in the previous inequality, using that $w_{i,j} = \xi_{i,j}(j+1)^{\mu}$, we derive that
\begin{equation}\label{Bterm2}
\Econd{B}{j} \leq \frac{\gamma_{\text{low}} G_{i,j}^2}{4 \Econd{w_{i,j}}{j}} + \frac{\kappa_g^2}{\varsigma\gamma_{\text{low}}} \Econd{\frac{g_{i,j}^2}{w_{i,j}^2}}{j}
\end{equation}

$\bullet$ Now consider the $C$ term. Again, Young's inequality with 
\[
x = \frac{|G_{i,j}|}{(j+1)^{\mu/2} \sqrt{\Econd{\xi_{i,j}}{j}}},
\ms
\lambda = \frac{\gamlow }{2}
\tim{and}
y =   \frac{|g_{i,j}| \Econd{|g_{i,j}|}{j} }{\sqrt{\Econd{\xi_{i,j}}{j}} \xi_{i,j} (j+1)^\mu}
\]
yields that
\begin{equation}\label{Cterm2}
	\Econd{C}{j} \leq \frac{\gamlow G_{i,j}^2}{4 \Econd{w_{i,j}}{j}} + \frac{\kappa_g^2}{\varsigma\gamlow} \Econd{\frac{g_{i,j}^2}{w_{i,j}^2}}{j},
\end{equation}
where we used that $|g_{i,j}| \leq \kappa_g$ and $\xi_{i,j}\geq \varsigma$. 
Taking $\Econd{.}{j}$ in \req{Abound}, using \req{Cterm2} and
\req{Bterm2} and injecting the obtained bound  into
\req{decreasbound}, we obtain  \req{avrdecrease} with \req{avrgi-kappas}.

\end{document}